\documentclass[a4paper,11pt,twoside,reqno]{amsart}

\usepackage[utf8]{inputenc}
\usepackage[plainpages=false,pdfpagelabels=true]{hyperref}
\usepackage{amssymb,amsthm}
\usepackage{slashed}
\usepackage{dsfont}

\newtheorem{Satz}{Theorem}[section]
\newtheorem{Prop}[Satz]{Proposition}
\newtheorem{Lem}[Satz]{Lemma}

\newtheorem{Cor}[Satz]{Corollary}
\newtheorem{Dfn}[Satz]{Definition}
\newtheorem{Bem}[Satz]{Remark}
\newtheorem{Bsp}[Satz]{Example}

\allowdisplaybreaks[1]

\renewcommand{\epsilon}{\varepsilon}
\newcommand{\cR}{{\mathcal R}}
\newcommand{\R}{\ensuremath{\mathbb{R}}}
\newcommand{\N}{\ensuremath{\mathbb{N}}}
\newcommand{\cG}{{\mathcal G}}
\newcommand{\D}{\slashed{D}}
\newcommand{\p}{\slashed{\partial}}
\newcommand{\sff}{\mathrm{I\!I}}
\newcommand{\hol}{\operatorname{Hol}}
\newcommand{\Hom}{\operatorname{Hom}}
\newcommand{\zarg}{(d\phi)^{\underline{2}}(\operatorname{vol}_h^\sharp)}

\numberwithin{equation}{section}

\title{Magnetic Dirac-harmonic Maps}
\author{Volker Branding}
\date{\today}
\address{TU Wien\\
Institut für diskrete Mathematik und Geometrie\\
Karlsplatz 13, 1040 Wien}
\email[]{volker@geometrie.tuwien.ac.at}
\urladdr{http://www.geometrie.tuwien.ac.at/branding/}
\begin{document}

\begin{abstract}
We study a functional, whose critical points couple Dirac-harmonic maps from surfaces with a two form. 
The critical points can be interpreted as coupling the prescribed mean curvature equation 
to spinor fields. On the other hand, this functional also arises as part of the supersymmetric sigma model in theoretical physics.
In two dimensions it is conformally invariant.
We call critical points of this functional magnetic Dirac-harmonic maps. 
We study geometric and analytic properties of magnetic Dirac-harmonic maps including their regularity
and the removal of isolated singularities.
\end{abstract} 

\maketitle

\section{Introduction}
Dirac-harmonic maps from Riemannian surfaces to a Riemannian manifold form 
a geometric variational problem with rich structure. They are critical points of an energy
functional that couples the equation for harmonic maps to spinor fields.
Thus, a Dirac-harmonic map is given by a pair \((\phi,\psi)\), where
\(\phi\colon M\to N\) is a map and \(\psi\) a vector spinor. 
For the case of a two-dimensional domain, Dirac-harmonic maps belong to the class of conformally invariant
variational problems. Due to their conformal invariance Dirac-harmonic maps
from surfaces share special properties and many substantial results have already been established.
They were introduced in \cite{MR2262709} together with a removable singularity theorem. 
The regularity of Dirac-harmonic maps from surfaces was first established for spherical targets
\cite{MR2176464}, which could be generalized to Dirac-harmonic maps to hypersurfaces \cite{MR2506243}.
The regularity in full generality was analysed in \cite{MR2544729} and extended recently in \cite{1306.4260}.

The energy identity for Dirac-harmonic maps from surfaces has been obtained in \cite{MR2267756}. 

Despite the fact that the analytical aspects of Dirac-harmonic maps are well understood at present,
the existence question is still not answered in full detail yet.
Using index-theoretical methods, Dirac-harmonic maps have been constructed in \cite{MR3070562}. These are uncoupled
in the sense that for a given harmonic map \(\phi_0\), a spinor \(\psi\) is constructed such that the pair \((\phi_0,\psi)\)
is a Dirac-harmonic map. A heat flow approach for Dirac-harmonic maps was recently studied in \cite{phd}.
An existence result for the boundary value problem for Dirac-harmonic maps was obtained in \cite{MR3044133}.

In this article we consider Dirac-harmonic maps coupled to a certain potential, which involves
a two-form on the target manifold \(N\). The potential we study is special in the sense that it does not break the conformal 
invariance. In the physical literature this potential can be interpreted as coupling the supersymmetric 
sigma model to a magnetic field. On the other hand in physics this extra part in the action is needed for the 
sake of anomaly cancellation. More on the physical background of the model we study can be found in \cite{MR1242500} and \cite{MR1910266}.
Our aim in this article is to show that this extension of Dirac-harmonic maps 
still has the many nice properties that Dirac-harmonic maps have.

Let us now describe the problem in more detail. For a map \(\phi\colon M\to N\)
we may study its differential \(d\phi\in\Gamma(T^*M\otimes\phi^{-1}TN)\), integrating the square of its
norm leads to the usual harmonic energy.

We assume that \((M,h)\) is a closed Riemannian spin surface with spinor bundle \(\Sigma M\), 
for more details about spin geometry see the book \cite{MR1031992}. Moreover, let \(N\) be another closed Riemannian
manifold of dimension \(n=\dim N\geq 3\). Together with the pullback bundle \(\phi^{-1}TN\) 
we can consider the twisted bundle \(\Sigma M\otimes\phi^{-1}TN\). 
The induced connection on this bundle will be denoted by \(\tilde{\nabla}\). 
Sections \(\psi\in\Gamma(\Sigma M\otimes\phi^{-1}TN)\) in this bundle  are called vector spinors and the natural operator acting on them 
is the twisted Dirac operator, denoted by \(\D\).  
It is an elliptic, first order operator, which is self-adjoint with respect to the \(L^2\)-norm.
More precisely, the twisted Dirac operator is given by \(\D=e_\alpha\cdot\tilde{\nabla}_{e_\alpha}\), where \(\{e_\alpha\}\) is
an orthonormal basis of \(TM\) and \(\cdot\) denotes Clifford multiplication. 
We are using the Einstein summation convention, that is we sum over repeated indices.
Clifford multiplication is skew-symmetric, namely
\[
\langle\chi,X\cdot\xi\rangle_{\Sigma M}=-\langle X\cdot\chi,\xi\rangle_{\Sigma M}
\]
for all \(\chi,\xi\in\Gamma(\Sigma M)\) and all \(X\in TM\).
Moreover, we use Greek letters for indices on \(M\) and Latin letters for indices on \(N\).
In terms of local coordinates the spinor \(\psi\) can be written as \(\psi=\psi^i\otimes\frac{\partial}{\partial y^i}\), 
and thus the twisted Dirac operator \(\D\) is locally given by
\[
\D\psi=\big(\p\psi^i+\Gamma^i_{jk}\nabla\phi^j\cdot\psi^k\big)\otimes\frac{\partial}{\partial y^i}.
\]
Here, \(\p\colon\Gamma(\Sigma M)\to\Gamma(\Sigma M)\) denotes the usual Dirac operator
and \(y^i\) are local coordinates on \(N\).

Let \(B\) be a two-form on the manifold \(N\), which we pullback by the map \(\phi\).
We now may state the central object of this article, which is the energy functional
\begin{equation}
\label{energy-functional}
E_B(\phi,\psi)=\frac{1}{2}\int_M(|d\phi|^2+\langle\psi,\D\psi\rangle+2\phi^{-1}B).
\end{equation}
Regarding the first term, the scalar product is taken with respect to the metric on the bundle
\(T^*M\otimes\phi^{-1}TN\), whereas in the second term we use the metric on \(\Sigma M\otimes\phi^{-1}TN\).

In local coordinates the energy functional acquires the form
\begin{equation}
E_B(\phi,\psi)=\frac{1}{2}\int_Mh^{\alpha\beta}g_{ij}\frac{\partial\phi^i}{\partial x_\alpha}\frac{\partial\phi^j}{\partial x_\beta}
+g_{ij}\langle\psi^i,(\D\psi)^j\rangle
+\epsilon^{\alpha\beta}B_{ij}\frac{\partial\phi^i}{\partial x_\alpha}\frac{\partial\phi^j}{\partial x_\beta}
, 
\end{equation}
where \(h_{\alpha\beta}\) denotes the Riemannian metric on \(M\) and \(g^{ij}\) the Riemannian metric on \(N\).
Moreover, \(\epsilon^{\alpha\beta}\) is the antisymmetric tensor in two dimensions and \(x_\alpha\) are coordinates on \(M\).

To state the critical points of the energy functional (\ref{energy-functional}), we need the following
definition. Let \(Z\in\Gamma(\operatorname{Hom}(\Lambda^2TN,TN))\cong\Gamma(\Lambda^2 T^*N\otimes TN)\)
be a tensor field satisfying
\begin{equation}
\Omega:=g(\cdot,Z(\cdot))
\end{equation}
with a three-form \(\Omega\).

The critical points of the energy functional (\ref{energy-functional}) are given by
\begin{eqnarray}
\label{euler-lagrange-phi-intro}
\tau(\phi)&=&\cR(\phi,\psi)+Z(\zarg), \\
\label{euler-lagrange-psi-intro}
\D\psi&=&0,
\end{eqnarray}
where \(\tau\) is the tension field of the map \(\phi\) and \(\zarg\) is defined in terms of an orthonormal basis \(\{e_1,e_2\}\) of \(M\) by
\(\zarg:=d\phi(e_1)\wedge d\phi(e_2)\). The quantity \(Z\) is the vector bundle homomorphism defined via the three-form \(\Omega=dB\) just above and
the curvature term \(\cR(\phi,\psi)\) is given by
\[
\cR(\phi,\psi)=\frac{1}{2}R^N(e_\alpha\cdot\psi,\psi)d\phi(e_\alpha),
\]
where \(R^N\) denotes the Riemann curvature tensor on \(N\).

We call solutions \((\phi,\psi)\) of (\ref{euler-lagrange-phi-intro}) and (\ref{euler-lagrange-psi-intro}) \emph{magnetic Dirac-harmonic maps}.
\par\medskip

This paper is organized as follows: In Section 2 we derive the Euler-Lagrange equations of (\ref{energy-functional})
and discuss their various limits. In Section 3 we introduce the energy-momentum tensor
of the energy functional (\ref{energy-functional}) and study the associated holomorphic
differential. Section 4 then establishes the regularity of magnetic Dirac harmonic maps 
and finally we prove the removable singularity theorem.

In order to prove our results, we make use of the analytical tools for Dirac-harmonic maps
provided in \cite{MR2262709} and \cite{MR2176464}. Since magnetic Dirac-harmonic maps have
the same analytic structure as Dirac-harmonic maps, many of the known results can be generalized easily.

\section{Magnetic Dirac-harmonic maps}
In this section we derive the critical points of the energy functional (\ref{energy-functional}),
analyze its limits and its connection to the literature. We will always assume that we have fixed
a spin structure on the surface \(M\).
\begin{Prop}[First Variation]
Let \((M,h)\) be a closed Riemannian spin surface and \((N,g)\) a closed Riemannian manifold.
Then the critical points of the energy functional (\ref{energy-functional}) are given by 
\begin{eqnarray}
\label{euler-lagrange-phi}
\tau(\phi)&=&\cR(\phi,\psi)+Z(\zarg), \\
\label{euler-lagrange-psi}
\D\psi&=&0
\end{eqnarray}
with \(\cR(\phi,\psi),Z(\zarg) \in\Gamma(\phi^{-1}TN)\) as defined in the introduction.
\end{Prop}

\begin{proof}
We choose a local orthonormal basis \(\{e_\alpha\}\) on \(M\) such that
\([e_\alpha,\partial_t]=0\) and also \(\nabla_{\partial_t}e_\alpha=0\) at a considered point.
We start by deriving the Euler-Lagrange equation for the spinor $\psi$. 
Therefore, we consider a variation of $\psi$ with $\phi$ fixed and $\frac{\tilde{\nabla}\psi_t}{\partial t}\big|_{t=0}=\chi$.
We find
\[
\frac{\delta}{\delta\psi}E_B(\phi,\psi_t)=
\frac{1}{2}\int_M(\langle\chi,\D\psi\rangle
+\langle\psi,\D\chi\rangle)dM
=\int_M\langle\chi,\D\psi\rangle dM.
\]
To derive the Euler-Lagrange equation for $\phi$, consider a family of smooth variations of $\phi$ satisfying 
$\frac{\partial\phi_t}{\partial t}\big|_{t=0}=\eta$, while keeping the $\psi^i$
in \(\psi(x)=\psi^i(x)\otimes\frac{\partial}{\partial y^i}(\phi_t(x))\) fixed.
The variation with respect to $\phi$ of the following terms has already
been computed in \cite{MR2262709}, p.\ 413, Prop. 2.1.:
\begin{eqnarray*}
\frac{\delta}{\delta\phi}\frac{1}{2}\int_M|d\phi_t|^2dM&=&-\int_M\langle\tau(\phi),\eta\rangle dM, \\
\frac{\delta}{\delta\phi}\frac{1}{2}\int_M\langle\psi,\D\psi\rangle dM&=&\int_M(\langle\frac{\nabla\psi}{\partial t},\D\psi\rangle
+\langle\cR(\phi,\psi),\eta\rangle)dM.
\end{eqnarray*}
The first variation of the action involving the two-form \(B\) leads to
\[
\frac{\delta}{\delta\phi}\int_M\phi_t^{-1}B=\int_M\langle\Omega(\eta,\zarg)\rangle dM,
\]
where \(\Omega=dB\) is a three-form on \(N\). A detailed derivation of this formula can be found in \cite{kohdiss}, Chapter 2.
As already explained in the introduction, we associate to a 3-form \(\Omega\in\Gamma(\Lambda^3T^*N)\) a smooth section
of \(\Hom(\Lambda^2TN,TN)\) via the vector bundle homomorphism \(Z:\Lambda^2TN\to TN\)
defined by the equation
\begin{equation}
\langle\eta,Z(\xi_1\wedge\xi_2)\rangle=\Omega(\eta,\xi_1,\xi_2)
\end{equation}
for all \(\eta,\xi_1,\xi_2\in TN\). Adding up the different contributions, we get
\begin{align*}
\frac{\delta}{\delta\phi}E_B(\phi_t,\psi)=\int_M\big(&\langle-\tau(\phi)+\cR(\phi,\psi)+Z(\zarg),\eta\rangle \\
&+\langle\D\psi,\frac{\nabla\psi}{\partial t}\big|_{t=0}\rangle\big)dM. 
\end{align*}

Using the Euler-Lagrange equation for $\psi$, which was deduced before, the result follows.
\end{proof}

\begin{Bem}
In principle, it is also possible to study the energy functional (\ref{energy-functional})
for higher dimensional \(M\), that is for \(m=\dim M\geq 3\). In this case, one would 
consider the pullback of an \(m\)-form on \(N\) in the energy functional and
the Euler-Lagrange equation for the map \(\phi\) would read
\[
\tau(\phi)=\cR(\phi,\psi)+Z((d\phi)^{\underline{m}}(\operatorname{vol}_h^\sharp)).
\]
For a derivation, see again \cite{kohdiss}, Chapter 2.
The non-linearity on the right hand side is even worse than the non-linearity in the harmonic map equation.
Thus, from an analytical point of view, one cannot expect to get a manageable problem. This is in contrast to Dirac-harmonic maps, which
form a reasonable problem in every dimension.
\end{Bem}

It is obvious that there exist several trivial solutions of (\ref{euler-lagrange-phi}) and (\ref{euler-lagrange-psi}).
This includes harmonic maps (for \(\psi=0\) and \(B=0\)), harmonic spinors (for \(\phi=0\))
and Dirac-harmonic maps (for \(B=0\)).

\begin{Bem}
The two form contribution in the energy functional 
can be formulated abstractly in terms of bundle gerbes and surfaces holonomy. 
Ignoring the spinor \(\psi\) for the moment, one studies the \(U(1)\)-valued functional
\begin{equation*}
e^{iE[\phi]}=e^{iE_h(\phi)}\hol_\cG(\phi)\text{ with } E_h(\phi)=\frac{1}{2}\int_M|d\phi|^2dM,
\end{equation*}
which is called Wess-Zumino term. For more details in this direction, see
\cite{MR2681604} and \cite{MR2648325}.
\end{Bem}

\begin{Bsp}
Assume that \(N=\R^3\) and denote \(2dB=H(y)dy_1dy_2dy_3\). Then the Euler-Lagrange equations
(\ref{euler-lagrange-phi}) and (\ref{euler-lagrange-psi}) acquire the form
\begin{eqnarray*}
\Delta\phi&=&2H(\phi)\partial_{x_1}\phi\wedge\partial_{x_2}\phi, \\
\p\psi^i&=&0, \qquad i=1,2,3.
\end{eqnarray*}
Hence, for the map \(\phi\) we obtain the so called \emph{H-surface system},
whereas the equation for \(\psi\) reduces to the equation for harmonic spinors on the surface \(M\).
\end{Bsp}
If the curvature of the target manifold \(N\) does not vanish and the map \(\phi\colon M\to N\) is an isometric embedding, 
then equations (\ref{euler-lagrange-phi}) and (\ref{euler-lagrange-psi}) naturally couple the prescribed mean
curvature equation to spinor fields.

\begin{Bsp}
If both \(M\) and \(N\) are Riemann surfaces, the magnetic term \(Z(\zarg)\) vanishes due to dimensional reasons.
In this case the Euler-Lagrange equations (\ref{euler-lagrange-phi}) and (\ref{euler-lagrange-psi})
reduce to the ones for Dirac-harmonic maps.
\end{Bsp}

\begin{Bem}
As already stated in the introduction, the study of the functional (\ref{energy-functional}) is motivated from
what physicists call the supersymmetric sigma model in two dimensions. 
In physics this model is usually formulated in terms of superfields. If one expands these into 
component fields one finds even more contributions that can be coupled to Dirac-harmonic maps.

One of these terms is given by
\[
\int_M\langle R^N(\phi)(\psi,\psi)\psi,\psi)\rangle dM
\]
and it also respects the conformal invariance. Dirac-harmonic maps coupled to this term (called \emph{Dirac-harmonic maps with curvature term}) are studied in \cite{MR2370260}.

In \cite{MR3018163} the authors couple Dirac-harmonic maps to a Ricci type potential
\[
\int_MR_{ij}(\phi)\langle\psi^i,\psi^j\rangle dM
\]
and establish the regularity of the critical points.
\end{Bem}

\section{Geometric aspects of magnetic Dirac-harmonic maps}
In this section we will study some geometric properties of magnetic Dirac-harmonic maps. 
First of all, we analyze the conformal invariance of the energy functional (\ref{energy-functional}).
\begin{Lem} 
The energy functional \(E_B(\phi,\psi)\) is conformally invariant. 
\end{Lem}
\begin{proof}
It is well known that all three contributions in (\ref{energy-functional}) are invariant under conformal transformations
of the metric on \(M\). For more details, see \cite{MR2262709}, Lemma 3.1 and \cite{MR744314}.
\end{proof}

For both harmonic maps/Dirac-harmonic maps from surfaces, there exists
a holomorphic quadratic differential, the so-called \emph{Hopf differential}.
We want to find an analogue for magnetic Dirac-harmonic maps. \\
Thus, let \((\phi,\psi)\) be a magnetic Dirac-harmonic map.
On a small domain \(\tilde{M}\) of \(M\) we choose a local isothermal
parameter \(z=x+iy\) and set
\begin{equation}
\label{hopf-differential}
T(z)dz^2=(|\phi_x|^2-|\phi_y|^2-2i\langle\phi_x,\phi_y\rangle+\langle\psi,\partial_x\cdot\tilde{\nabla}_{\partial_x}\psi\rangle-i\langle\psi,\partial_x\cdot\tilde{\nabla}_{\partial_y}\psi\rangle)dz^2
\end{equation}
with \(\partial_x=\frac{\partial}{\partial x}\) and \(\partial_y=\frac{\partial}{\partial y}\).

\begin{Dfn}
We define a two-tensor by 
\begin{equation}
\label{energy-momentum-tensor}
T_{\alpha\beta}:=2\langle d\phi(e_\alpha),d\phi(e_\beta)\rangle-\delta_{\alpha\beta}|d\phi|^2+\langle\psi,e_\alpha\cdot\tilde{\nabla}_{e_\beta}\psi\rangle.
\end{equation}
The tensor \(T_{\alpha\beta}\) is called the \emph{energy-momentum tensor}.
\end{Dfn}
It is easy to see that \(T_{\alpha\beta}\) is symmetric and traceless, when \((\phi,\psi)\) is a magnetic Dirac-harmonic map.
\begin{Bem}
Note that the energy-momentum tensor \(T_{\alpha\beta}\) for magnetic Dirac-harmonic maps is the same as the energy-momentum tensor for Dirac-harmonic maps.
This is not surprising since the energy-momentum tensor is calculated by varying the action functional with respect to the metric \(h\) on \(M\).
Since the pullback of the two-form \(B\) does not depend on the metric on \(M\) we do not get a new contribution. 
\end{Bem}

\begin{Lem}
Let \((\phi,\psi)\in C^2(M,N)\times C^2(M,\Sigma M\otimes\phi^{-1}TN)\) be a solution of (\ref{euler-lagrange-phi}) and (\ref{euler-lagrange-psi}).
Then the energy-momentum tensor (\ref{energy-momentum-tensor}) is covariantly conserved, that is
\[
\nabla_{e_\alpha}T_{\alpha\beta}=0.
\]
\end{Lem}
\begin{proof}
We choose a local orthonormal basis \(\{e_\alpha\}\) on \(M\) such that
\([e_\alpha,e_\beta]=0\) and also \(\nabla_{e_\alpha}e_\beta=0\) at a considered point.
Set
\[
C_{\alpha\beta}:=2\langle d\phi(e_\alpha),d\phi(e_\beta)\rangle-\delta_{\alpha\beta}|d\phi|^2, \qquad D_{\alpha\beta}:=\langle\psi,e_\alpha\cdot\tilde{\nabla}_{e_\beta}\psi\rangle.
\]
First, we compute
\begin{eqnarray*}
\nabla_{e_\alpha}C_{\alpha\beta}&=&2\langle\tau(\phi),d\phi(e_\beta)\rangle\\
&=&2\langle\cR(\phi,\psi),d\phi(e_\beta)\rangle+\underbrace{\langle Z(\zarg),d\phi(e_\beta)\rangle}_{=\Omega(d\phi(e_1),d\phi(e_2),d\phi(e_\beta))=0}\\
&=&2\langle\cR(\phi,\psi),d\phi(e_\beta)\rangle.
\end{eqnarray*}
On the other hand, we find
\begin{eqnarray*}
\nabla_{e_\alpha}D_{\alpha\beta}&=&\nabla_{e_\alpha}\langle\psi,e_\alpha\cdot\tilde{\nabla}_{e_\beta}\psi\rangle \\
&=&\langle\tilde{\nabla}_{e_\alpha}\psi,e_\alpha\cdot\tilde{\nabla}_{e_\beta}\psi\rangle +\langle\psi,\D\tilde{\nabla}_{e_\beta}\psi\rangle \\
&=&-\langle\D\psi,\tilde{\nabla}_{e_\beta}\psi\rangle +\langle\psi,\D\tilde{\nabla}_{e_\beta}\psi\rangle \\
&=&\langle\psi,\D\tilde{\nabla}_{e_\beta}\psi\rangle,
\end{eqnarray*}
where we used that \(\psi\) is a solution of (\ref{euler-lagrange-psi}). 
The curvature tensor on the twisted bundle \(\Sigma M\otimes\phi^{-1}TN\) can be decomposed as
\[
\tilde{R}(e_\alpha,e_\beta)=R^{\Sigma M}(e_\alpha,e_\beta)\otimes\mathds{1}_{\phi^{-1}TN}
+\mathds{1}_{\Sigma M}\otimes R^N(d\phi(e_\alpha),d\phi(e_\beta)).
\]
Then, again, a direct calculation yields
\begin{eqnarray*}
\langle\psi,\D\tilde{\nabla}_{e_\beta}\psi\rangle&=&\langle\psi,\tilde{\nabla}_{e_\beta}\underbrace{\D\psi}_{=0}\rangle
+\underbrace{\langle\psi,e_\alpha\cdot R^{\Sigma M}(e_\alpha,e_\beta)\psi\rangle}_{=\frac{1}{2}\langle\psi,\operatorname{Ric}(e_\beta)\cdot\psi\rangle=0} \\
&&+\langle\psi,e_\alpha\cdot R^{N}(d\phi(e_\alpha),d\phi(e_\beta))\psi\rangle \\
&=&-2\langle\cR(\phi,\psi),d\phi(e_\beta)\rangle,
\end{eqnarray*}
which completes the proof.
\end{proof}

\begin{Prop}
\label{prop-energy-momentum-preserved}
The quadratic differential \(T(z)dz^2\) is holomorphic.
\end{Prop}
\begin{proof}
This follows directly from the last Lemma.
\end{proof}
For the rest of this section we assume that the spinor \(\psi\) vanishes such that we can draw a comparison with
so-called magnetic geodesics, see for example \cite{MR2788659}. These are the lower-dimensional analogue of magnetic Dirac-harmonic maps with vanishing spinors.
It is known that magnetic geodesics have constant energy. In the two-dimensional case we have a similar result,
but here we have to make curvature assumptions.

\begin{Lem}
Let \((M,h)\) be a closed Riemann surface with \(\operatorname{Ric}^M\geq 0\)
and suppose \(\psi=0\). Moreover, assume that the sectional curvature \(K^N\) of the target manifold \(N\) is non-positive,
then a solution \(\phi\in C^2(M,N)\) of (\ref{euler-lagrange-phi}) has constant energy.
\end{Lem}
\begin{proof}
Since \(\psi=0\), the map \(\phi\) satisifes \(\tau(\phi)=Z(\zarg).\)
Moreover, we note that
\begin{align*}
\langle\nabla_{e_\alpha}\tau(\phi),d\phi(e_\alpha)\rangle
&=\langle\nabla_{e_\alpha}Z(\zarg),d\phi(e_\alpha)\rangle \\
&=-\langle Z(\zarg),\tau(\phi)\rangle
=-|\tau(\phi)|^2.
\end{align*}
As a second step, we calculate
\begin{eqnarray*}
\Delta\frac{1}{2}|d\phi|^2&=&|\nabla d\phi|^2
-\langle R^N(d\phi(e_\alpha),d\phi(e_\beta))d\phi(e_\alpha),d\phi(e_\beta)\rangle\\
&&+\langle d\phi(\operatorname{Ric}^M(e_\beta)),d\phi(e_\beta)\rangle+\langle\nabla_{e_\alpha}\tau(\phi),d\phi(e_\alpha)\rangle\\
&\geq& |\nabla d\phi|^2 -|\tau(\phi)|^2\\
&\geq& 0.
\end{eqnarray*}
The statement then follows by application of the maximum principle.
\end{proof}
\begin{Cor}
Under the assumptions of the last Lemma, we also know that any solution \(\phi\in C^2(M,N)\) of (\ref{euler-lagrange-phi})
is totally geodesic. This follows from the fact that \(|d\phi|^2\) is constant and thus \(|\nabla d\phi|=0\).
\end{Cor}

For the further analysis it turns out to be useful to apply the Nash embedding theorem to
isometrically embed the manifold \(N\) in some \(\R^q\) of sufficiently high dimension.
The Euler-Lagrange equations for magnetic Dirac-harmonic maps then acquire the form
\begin{eqnarray}
\label{euler-lagrange-u}
-\Delta\phi&=&\sff(d\phi,d\phi)+Z(\zarg)+P(\sff(d\phi(e_\alpha),e_\alpha\cdot\psi),\psi),\\
\label{euler-lagrange-psirq}
\p\psi&=&\sff(e_\alpha\cdot\psi,d\phi(e_\alpha)),
\end{eqnarray}
where \(\sff\) denotes the second fundamental form of \(\phi\) in \(\R^q\) and \(P\) the shape operator.
For a detailed derivation see \cite{MR2176464} and \cite{kohdiss}.
Now, we have that \(\phi\colon M\to\R^q\) with \(\phi(x)\in N\). The vector spinor \(\psi\)  becomes
a vector of untwisted spinors \(\psi^1,\psi^2,\ldots,\psi^q\), more precisely \(\psi\in\Gamma(\Sigma M\otimes T\R^q)\).
The condition that \(\psi\) is along the map \(\phi\) is now encoded as
\[
\sum_{i=1}^q\nu^i\psi^i=0\qquad \text{for any normal vector }\nu \text{ at } \phi(x).
\]
Moreover, the vector bundle homomorphism \(Z\) can be extended to \(\R^q\)
by parallel transport. We may now turn to the study of the analytic aspects of magnetic Dirac-harmonic maps.

\section{Analytic aspects of magnetic Dirac-harmonic maps}
\subsection{Regularity of magnetic Dirac-harmonic maps}\mbox{}\\
In this section we want to study the regularity of magnetic Dirac-harmonic maps. 
The regularity question for the various limits of the energy functional \(E_B(\phi,\psi)\) is already fully developed.
For Dirac-harmonic maps, see \cite{MR2176464}, \cite{MR2506243} and \cite{MR2544729}. On the other hand, for \(\psi=0\) the issue
of regularity has been analyzed extensively, see \cite{MR636440}, \cite{MR1331835}, \cite{MR1168525} and \cite{MR1913803} p.187 ff.

First of all, we need the definition of a weak solution of (\ref{euler-lagrange-phi}) and (\ref{euler-lagrange-psi}). 
Therefore, we define
\begin{eqnarray*}
\chi(M,N)&:=&\{(\phi,\psi)\in W^{1,2}(M,N)\times W^{1,\frac{4}{3}}(M,\Sigma M\otimes\phi^{-1}TN) \\
&&\hspace{0.5cm}\text{ with } (\ref{euler-lagrange-u}) \text{ and } (\ref{euler-lagrange-psirq}) \text{ a.e.}\}.
\end{eqnarray*}
\begin{Dfn}[Weak magnetic Dirac-harmonic Map]
A pair \((\phi,\psi)\in\chi(M,N)\) is called \emph{weak magnetic Dirac-harmonic map} from \(M\) to \(N\) if and only if the pair \((\phi,\psi)\) satisfies
\begin{eqnarray*}
&&\int_M(\langle\nabla\phi,\nabla\xi\rangle-\langle\sff(d\phi,d\phi),\xi\rangle+\langle P(\sff(d\phi(e_\alpha),e_\alpha\cdot\psi),\psi),\xi\rangle \\
&&\hspace{0.5cm}+\langle Z(\zarg),\xi\rangle)dM=0, \\
&&\int_M(\langle\p\psi,\eta\rangle-\langle\sff(e_\alpha\cdot\psi,d\phi(e_\alpha)),\eta\rangle)dM=0
\end{eqnarray*}
for all \(\xi\in C^\infty(M,\R^q)\) and \(\eta\in C^\infty(M,\Sigma M\otimes T\R^q)\).
\end{Dfn}

In the upcoming analysis we will make use of the following powerful result due to Rivi\`ere (see \cite{MR2285745}):
\begin{Satz}
\label{theorem-riviere}
Let \(D\) be the unit disc in \(\R^2\) and \(q\in\N\) a fixed number. Then for every \(A={A_{~j}^i},1\leq i,j\leq q\)
in \(L^2(D,so(q)\otimes\R^2)\)(that is for all \(i,j\in 1,\ldots q, A_{~j}^i\in L^2(D,\R^2)\)
and \(A_{~j}^i=-A_{~i}^j\)), every \(\phi\in W^{1,2}(D,\R^q)\) solving
\begin{equation}
-\Delta\phi=A\cdot\nabla\phi
\end{equation}
is continuous. 
The notation should be understood as \(-\Delta \phi^i=\sum_{j=1}^qA_{~j}^i\cdot\nabla\phi^j\) for all \(1\leq i\leq q\).
\end{Satz}
To apply the above Theorem to magnetic Dirac-harmonic maps, 
we have to rewrite the Euler-Lagrange equations. 
Thus, let us fix the notation, we follow the presentation in \cite{MR2506243} for Dirac-harmonic maps.
We need to express all data in terms of the ambient space \(\R^q\).
We denote coordinates in the ambient space \(\R^q\) by \((y^1,y^2,\ldots,y^{q})\).
Let \(\nu_l,l=n+1,\ldots,q\) be an orthonormal frame field for the normal bundle \(T^\perp N\).
For \(X,Y\in T_yN\) and \(\nabla_Y\nu_k=Y^i\frac{\partial\nu_k}{\partial y^i}\)
we have
\begin{equation}
\sff_y(X,Y)=\langle X,\nabla_Y\nu_l\rangle\nu_l=X^iY^j\frac{\partial\nu^i_l}{\partial y^j}\nu_l.
\end{equation}
Let \(D\) be a domain in \(M\) and consider a weak magnetic Dirac-harmonic map \((\phi,\psi)\in\chi(M,N)\).
We choose local isothermal coordinates \(z=x+iy\), set \(e_1=\partial_x\), \(e_2=\partial_y\)
and use the notation \(\phi_\alpha=d\phi(e_\alpha)\). 
Moreover, note that \(\phi_\alpha\in TN\) and \(\nu_l\in T^\perp N\), which implies that
\begin{equation}
\label{orthogonality}
\phi^i_\alpha\nu^i_l=0
\end{equation}
for all \(\alpha\). Hence, we may write
\begin{equation}
\label{skew-sff}
\sff^m(\phi_\alpha,\phi_\alpha)=\phi^i_\alpha\phi^j_\alpha\left(\frac{\partial\nu_l^i}{\partial y^j}\nu_l^m-\frac{\partial\nu_l^m}{\partial y^j}\nu_l^i\right),\qquad m=1,2,\ldots,q, 
\end{equation}
where we used (\ref{orthogonality}) in the second term on the right hand side.
Following \cite{BVP}, p.7, the term on the right hand side of (\ref{euler-lagrange-u}) involving the shape operator can also
be written in a skew-symmetric way, namely
\begin{align}
\label{skew-p}
\operatorname{Re} P^m(\sff(\phi_\alpha,e_\alpha\cdot\psi),\psi)=& \\
\nonumber\phi^i_\alpha\langle\psi^k,e_\alpha\cdot\psi^j\rangle &
\Bigg(\bigg(\frac{\partial\nu_l}{\partial y^j}\bigg)^{\top,i}\bigg(\frac{\partial\nu_l}{\partial y^k}\bigg)^{\top,m}-\bigg(\frac{\partial\nu_l}{\partial y^k}\bigg)^{\top,i}\bigg(\frac{\partial\nu_l}{\partial y^j}\bigg)^{\top,m}\Bigg).
\end{align}
Here, \(\top\) denotes the projection map \(\top\colon\R^q\to T_yN\).
Finally, we note that
\[
Z^m(\zarg)=Z^m(\partial_{y^i}\wedge\partial_{y^j})\phi^i_x\phi^j_y.
\]
By the definition of \(Z\) and exploiting the skew-symmetry of the three-form \(\Omega\), we find (see also \cite{MR2285745}, p.13)
\begin{equation}
\label{skew-z}
Z^k(\partial_{y^i}\wedge\partial_{y^j})=-Z^i(\partial_{y^k}\wedge\partial_{y^j}).
\end{equation}
We are now in the position to show that magnetic Dirac-harmonic maps have a structure such that Theorem \ref{theorem-riviere}
can be applied. Throughout the upcoming calculation we will assume that
both \(|Z|_{L^\infty}\leq C\) and \(|\nabla Z|_{L^\infty}\leq C\).
\begin{Prop}
Let \((M,h)\) be a closed Riemannian spin surface and let \(N\) be a compact Riemannian manifold. Assume that \((\phi,\psi)\in\chi(M,N)\) is a weak solution of (\ref{euler-lagrange-u})
and (\ref{euler-lagrange-psirq}). Let \(D\) be a simply connected domain of \(M\). Then there exists \(A^i_{~j}\in L^2(D,so(q)\otimes\R^2)\)
such that
\begin{equation}
-\Delta\phi^m=A^m_{~i}\cdot\nabla\phi^i
\end{equation}
holds.
\end{Prop}
\begin{proof}
By assumption \(N\subset\R^{q}\) is compact, we denote its unit normal field by \(\nu_l,l=n+1,\ldots,q\).
Exploiting the skew-symmetry of (\ref{skew-sff}), (\ref{skew-p}) and (\ref{skew-z}), we denote
\[
A^m_{~i}=\begin{pmatrix}
f^m_{~i} \\
g^m_{~i}
\end{pmatrix},
\qquad
i,m=1,2,\ldots,q
\]
with
\begin{align*}
f^m_{~i}:=&\left(\frac{\partial\nu_l^i}{\partial y^j}\nu_l^m-\frac{\partial\nu_l^m}{\partial y^j}\nu_l^i\right)\phi^j_x+Z^m(\partial_{y^i}\wedge\partial_{y^j})\phi^j_y\\
&+\langle\psi^k,\partial_x\cdot\psi^j\rangle
\Bigg(\bigg(\frac{\partial\nu_l}{\partial y^j}\bigg)^{\top,i}\bigg(\frac{\partial\nu_l}{\partial y^k}\bigg)^{\top,m}-\bigg(\frac{\partial\nu_l}{\partial y^k}\bigg)^{\top,i}\bigg(\frac{\partial\nu_l}{\partial y^j}\bigg)^{\top,m}\Bigg),\\
g^m_{~i}:=&\left(\frac{\partial\nu_l^i}{\partial y^j}\nu_l^m-\frac{\partial\nu_l^m}{\partial y^j}\nu_l^i\right)\phi^j_y-Z^m(\partial_{y^i}\wedge\partial_{y^j})\phi^j_x\\
&+\langle\psi^k,\partial_y\cdot\psi^j\rangle
\Bigg(\bigg(\frac{\partial\nu_l}{\partial y^j}\bigg)^{\top,i}\bigg(\frac{\partial\nu_l}{\partial y^k}\bigg)^{\top,m}-\bigg(\frac{\partial\nu_l}{\partial y^k}\bigg)^{\top,i}\bigg(\frac{\partial\nu_l}{\partial y^j}\bigg)^{\top,m}\Bigg).
\end{align*}
Now, we can write (\ref{euler-lagrange-u}) in the following form
\[
-\Delta\phi^m=A^m_{~i}\cdot\nabla\phi^i.
\]
It remains to show that \(A^m_{~i}\in L^2(D,so(q)\otimes\R^2)\).
This follows directly since the pair \((\phi,\psi)\) is a weak solution of (\ref{euler-lagrange-u}), (\ref{euler-lagrange-psirq}) 
and the Sobolev embedding \(|\psi|_ {L^4}\leq C|\psi|_{W^{1,\frac{4}{3}}}\).
The skew-symmetry of \(A^m_{~i}\) can be read of from its definition and the properties of \(Z\), see (\ref{skew-z}).
\end{proof}

\begin{Cor}
Let \((M,h)\) be a closed Riemannian spin surface and let \(N\subset \R^{q}\) be a compact manifold.
Suppose that \((\phi,\psi)\in\chi(M,N)\) is a weak magnetic Dirac-harmonic map. 
Then by the last Proposition and Theorem \ref{theorem-riviere}, 
we may deduce that \(\phi^m\) is continuous, \(m=1,2,\ldots,q\), hence \(\phi\in C^0(M,N)\).
\end{Cor}

We are now in the position to apply the regularity theory developed for Dirac-harmonic maps in \cite{MR2176464}.
Since \(\phi\) is  in \(C^0(D,N)\) we may choose local coordinates \(\{y_i\}\) on \(N\).
In these coordinates equations (\ref{euler-lagrange-u}) and (\ref{euler-lagrange-psirq}) acquire the form
\begin{eqnarray}
\label{euler-lagrange-u-coordinates}
\Delta\phi^m&=&-\Gamma^m_{ij}(\phi)\phi^i_\alpha\phi^j_\alpha+\frac{1}{2}R^m_{~kjl}(\phi)\langle\psi^k,\nabla\phi^j\cdot\psi^l\rangle \\
&&\nonumber+Z^m(\phi)(\partial_{y^i}\wedge\partial_{y^j})\phi^i_x\phi^j_y,\\
\label{euler-lagrange-psi-coordinates}
\p\psi^m&=&-\Gamma^m_{ij}(\phi)\nabla\phi^j\cdot\psi^i.
\end{eqnarray}
As a next step we need two establish two auxiliary Lemmas. 
These are very similar to Lemma 2.4 and Lemma 2.5 in \cite{MR2176464}.
Hence, we will not prove both in full detail. 

\begin{Lem}
\label{lemma1-regularity}
Let the pair \((\phi,\psi)\) be a weak solution of (\ref{euler-lagrange-u-coordinates}) and (\ref{euler-lagrange-psi-coordinates}).
If \(\phi\in C^0\cap W^{1,2}(D,N)\), then for any \(\epsilon>0\), there is a \(\rho>0\) such that
\begin{equation}
\int_{D(x_1,\rho)}|\nabla\phi|^2\eta^2\leq\epsilon\int_{D(x_1,\rho)}|\nabla\eta|^2+C\epsilon\bigg(\int_{D(x_1,\rho)}|\psi|^4\eta^4\bigg)^\frac{1}{2},
\end{equation}
\end{Lem}
where \(D(x_1,\rho)\subset D, \eta\in W^{1,2}_0(D(x_1,\rho),\R)\) and \(C\) is a positive constant
independent of \(\epsilon,\rho, \phi\) and \(\psi\).
\begin{proof}
Since \(\phi\in C^0(D,N)\) we may choose local coordinates and set
\[
G(x,\phi,d\phi,\psi):=\Gamma(\phi)(d\phi,d\phi)-Z(\phi)(\zarg)-\frac{1}{2}R^N(\phi)(e_\alpha\cdot\psi,\psi)d\phi(e_\alpha).
\]
Then the weak form of (\ref{euler-lagrange-u-coordinates}) is
\begin{equation}
\label{lemma1-regularity-a}
\int_{D(x_1,\rho)}\langle\nabla\phi,\nabla\xi\rangle=\int_{D(x_1,\rho)}\langle G(x,\phi,d\phi,\psi),\xi\rangle
\end{equation}
for any \(\xi\in W^{1,2}\cap L^\infty(D,\R^q)\). Now choose \(\xi(x)=(\phi(x)-\phi(x_1))\eta^2\), then we get
\begin{equation}
\label{lemma1-regularity-b}
\int_{D(x_1,\rho)}\langle\nabla\phi,\nabla\xi\rangle=\int_{D(x_1,\rho)}|d\phi|^2\eta^2+2\int_{D(x_1,\rho)}\eta(\phi(x)-\phi(x_1))\langle\nabla\phi,\nabla\eta\rangle.
\end{equation}
On the other hand, we find
\begin{align}
\label{lemma1-regularity-c}
\int_{D(x_1,\rho)}&\langle G(x,\phi,d\phi,\psi),\xi\rangle \\
=&\nonumber\int_{D(x_1,\rho)}\langle\Gamma(\phi)(d\phi,d\phi)-Z(\phi)(\zarg)\\
\nonumber &\hspace{1.5cm}-\frac{1}{2}R^N(\phi)(e_\alpha\cdot\psi,\psi)d\phi(e_\alpha),\eta^2(\phi(x)-\phi(x_1))\rangle\\
\leq&\nonumber C\epsilon_1\sup_{D(x_1,\rho)}|\phi(x)-\phi(x_1)|\int_{D(x_1,\rho)}|d\phi|^2\eta^2\\
&\nonumber+C\sup_{D(x_1,\rho)}|\phi(x)-\phi(x_1)|\bigg(\int_{D(x_1,\rho)}|d\phi|^2\bigg)^\frac{1}{2}\bigg(\int_{D(x_1,\rho)}|\psi|^4\eta^4\bigg)^\frac{1}{2},
\end{align}
where \(\epsilon_1>0\) is a given small number. We omitted some details in the calculation, see the proof
of Lemma 2.4 in \cite{MR2176464}.
In addition, we have
\begin{align}
\label{lemma1-regularity-d}
2\int_{D(x_1,\rho)}&\eta(\phi(x)-\phi(x_1))\langle\nabla\phi,\nabla\eta\rangle\\
\nonumber\leq&C\sup_{D(x_1,\rho)}|\phi(x)-\phi(x_1)|\int_{D(x_1,\rho)}|d\phi||\nabla\eta|\eta\\
\nonumber\leq&\frac{1}{2}\int_{D(x_1,\rho)}|d\phi|^2\eta^2+8\sup_{D(x_1,\rho)}|\phi(x)-\phi(x_1)|^2\int_{D(x_1,\rho)}|\nabla\eta|^2.
\end{align}
Substituting (\ref{lemma1-regularity-b}), (\ref{lemma1-regularity-c}), (\ref{lemma1-regularity-d}) into (\ref{lemma1-regularity-a}) and choosing \(\rho\) small enough then proves the Lemma.
\end{proof}

\begin{Lem}
\label{lemma2-regularity}
Let \(\phi\in C^0\cap W^{1,4}(D(x_0,R),N)\) and suppose that \((\phi,\psi)\) is a weak solution of (\ref{euler-lagrange-u-coordinates}) and (\ref{euler-lagrange-psi-coordinates}). 
Then for \(R\) sufficiently small, we have
\begin{equation}
\label{lemma2-regularity-result}
|\nabla^2\phi|_{L^2(D(x_0),\frac{R}{2})}+|d\phi|^2_{L^4((D(x_0),\frac{R}{2})}\leq C_1|d\phi|^2_{L^2(D(x_0,R))},
\end{equation}
where \(C_1\) is a positive constant depending on \(|\phi|_{C^0(D,N)}\) and \(R\).
\end{Lem}

\begin{proof}
We denote \(B:=D(x_0,R)\). For any \(\zeta\in W_0^{1,2}(B,\R^q)\) we have
\begin{equation}
\label{lemma2-regularity-a}
\int_B\langle\nabla\phi,\nabla\zeta\rangle=-\int_B\langle\Delta\phi,\zeta\rangle=\int_B\langle G(x,\phi,d\phi,\psi),\zeta\rangle.
\end{equation}
Choosing \(\zeta=\nabla_\gamma(\xi^2\nabla_\gamma\phi)\), where \(\xi\in C^\infty\cap W_0^{1,2}(B,\R)\) is to be determined later, a direct calculation yields
\begin{align}
\label{lemma2-regularity-b}
\int_B\langle\nabla_\gamma(\nabla_\beta\phi),\nabla_\beta(\xi^2\nabla_\gamma\phi)\rangle=\int_B\langle\nabla G,\nabla\phi\rangle\xi^2.
\end{align}
Moreover, note that
\begin{align}
\label{lemma2-regularity-c}
\nonumber\langle\nabla_\gamma(\nabla_\beta\phi),\nabla_\beta(\xi^2\nabla_\gamma\phi)\rangle
&=|\nabla_\gamma\nabla_\beta\phi|^2\xi^2+\langle\nabla_\gamma\nabla_\beta\phi,\nabla_\gamma\phi\rangle\nabla_\beta\xi^2 \\
&\geq|\nabla^2\phi|\xi^2-2|\nabla^2\phi||d\phi||\xi\nabla\xi|
\end{align}
and
\begin{eqnarray}
\label{lemma2-regularity-d}
|\langle\nabla G,\nabla\phi\rangle|&\leq &C(|d\phi|^4+|\nabla^2\phi||d\phi|^2+|\psi||\nabla\psi||d\phi|^2\\
&&\nonumber +|d\phi|^3|\psi|^2+|\nabla^2\phi||d\phi||\psi|^2).
\end{eqnarray}
Substituting (\ref{lemma2-regularity-c}) and (\ref{lemma2-regularity-d}) into (\ref{lemma2-regularity-b}) leads to the following inequality
\begin{align}
\nonumber \int_B|\nabla^2\phi|\xi^2\leq &C\bigg(\int_B|\nabla^2\phi||d\phi||\xi\nabla\xi|+\int_B|\nabla^2\phi||d\phi|^2\xi^2+\int_B|d\phi|^4\xi^2\\
\nonumber&+\int_B|\psi||\nabla\psi||d\phi|^2\xi^2+\int_B|d\phi|^3|\psi|^2\xi^2+\int_B|\nabla^2\phi||d\phi||\psi|^2\xi^2\bigg) \\
 &:= I+II+III+IV+V+VI.
\end{align}
The first two terms \(I\) and \(II\) can easily be estimated as
\begin{align*}
C\int_B|\nabla^2\phi||d\phi||\xi\nabla\xi|&\leq C\delta_1\int_B|\nabla^2\phi|\xi^2+\frac{C}{\delta_1}\int_B|d\phi|^2|\nabla\xi|^2, \\
C\int_B|\nabla^2\phi||d\phi|^2\xi^2 &\leq C\delta_1\int_B|\nabla^2\phi|\xi^2 +\frac{C}{\delta_1}\int_B|d\phi|^4\xi^2,
\end{align*}
where \(\delta_1>0\) is a small constant. The term \(III\) already has the shape that we need.
The last three contributions \(IV,V,VI\) can be estimated the same way as in the proof of Lemma 2.5 in \cite{MR2176464} leading to
\begin{align*}
IV&\leq \frac{1}{8}\int_B|\nabla^2\phi|\xi^2+\frac{1}{8}\int_B|d\phi|^2|\nabla\xi|^2, \\
V&\leq \frac{1}{8}\int_B|\nabla^2\phi|\xi^2+\frac{1}{8}\int_B|d\phi|^2|\nabla\xi|^2+C\int_B|d\phi|^4\xi^2, \\
VI&\leq \frac{1}{8}\int_B|\nabla^2\phi|\xi^2+\frac{1}{8}\int_B|d\phi|^2|\nabla\xi|^2,
\end{align*}
where we applied a local estimate for the spinor \(\psi\) derived in \cite{MR2176464}.
Applying these estimates, we find
\begin{equation}
\label{lemma2-regularity-e}
\int_{D(x_0,R)}|\nabla^2\phi|^2\xi^2\leq C\bigg(\int_{D(x_0,R)}|d\phi|^2|\nabla\xi|^2+\int_{D(x_0,R)}|d\phi|^4\xi^2\bigg).
\end{equation}
Now, for \(\epsilon>0\), let \(\rho>0\) be as in Lemma \ref{lemma1-regularity}.
Suppose that \(D(x_1,\rho)\subset D(x_0,R)\) and choose a cut-off function
\(\xi\in C^\infty_0(D(x_1,\rho)), 0\leq\xi\leq 1\) such that
\[
\xi=1 \text{ in } D(x_1,\frac{\rho}{2}),\qquad |\nabla\xi|\leq\frac{C}{\rho}\text{ in } D(x_1,\rho).
\]
For simplicity, we denote \(D_\rho:=D(x_1,\rho)\) and calculate
\begin{align}
\label{lemma2-regularity-f}
\nonumber\int_{D_\rho}|d\phi|^4\xi^2&=\int_{D_\rho}|d\phi|^2(|d\phi|\xi)^2\\
\nonumber&\leq\epsilon\int_{D_\rho}|\nabla(|d\phi|\xi)|^2 +C\epsilon\bigg(\int_{D_\rho}|\psi|^4|d\phi|^4\xi^4\bigg)^\frac{1}{2} \qquad \text{by Lemma \ref{lemma1-regularity}}\\
&\leq\epsilon\int_{D_\rho}|\nabla^2\phi|^2\xi^2+\epsilon\int_{D_\rho}|d\phi|^2|\nabla\xi|^2+C\epsilon\bigg(\int_{D_\rho}|\psi|^4|d\phi|^4\xi^4\bigg)^\frac{1}{2}.
\end{align}
As a next step we state an inequality for solutions of (\ref{euler-lagrange-psi-coordinates}), which was derived in the proof of Lemma 2.5 in \cite{MR2176464}
\[
\bigg(\int_{D_\rho}|\psi|^4|d\phi|^4\xi^4\bigg)^\frac{1}{2}\leq C\bigg(\int_{D_\rho}|\psi|^4\bigg)^\frac{1}{2}\bigg(\int_{D_\rho}|\nabla^2\phi|^2\xi^2+\int_{D_\rho}|d\phi|^2|\nabla\xi|^2\bigg).
\]
Using this estimate and choosing \(\rho\) such that \(|\psi|^2_{L^4(D_\rho)}\) is small enough, combining (\ref{lemma2-regularity-e}) and (\ref{lemma2-regularity-f}), we find
\[
\int_{D_\rho}|\nabla^2\phi|^2\xi^2\leq C\int_{D_\rho}|d\phi|^2|\nabla\xi|^2
\]
and thus
\begin{equation}
\label{lemma2-regularity-g}
\int_{D(x_1,\frac{\rho}{2})}|\nabla^2\phi|^2\leq\frac{C}{\rho^2}\int_{D(x_1,\rho)}|d\phi|^2.
\end{equation}
Covering \(D(x_0,\frac{R}{2})\) with \(\{D(x_1,\frac{\rho}{2})\}\) together with (\ref{lemma2-regularity-g}) gives the result.
\end{proof}

As for Dirac-harmonic maps \cite{MR2176464}, we may thus follow:

\begin{Satz}
Let \((M,h)\) be a closed Riemann spin surface. Suppose that the pair \((\phi,\psi)\colon(D,\delta_{\alpha\beta})\to(N,g_{ij})\) is a weak magnetic
Dirac-harmonic map. If \(\phi\) is continuous, then the pair \((\phi,\psi)\) is smooth.
\end{Satz}
\begin{proof}
First of all, we note that \(\phi\in W^{2,2}\cap W^{1,4}(D(x_0,\frac{R}{2}),N)\).
This can be achieved by replacing weak derivatives by difference quotients in the proof of Lemma \ref{lemma1-regularity}.

We may now establish the regularity of solutions of (\ref{euler-lagrange-u-coordinates}) and (\ref{euler-lagrange-psi-coordinates}).
Since \(\phi\in W^{2,2}\), we have by the Sobolev embedding theorem that \(\phi\in W^{1,p}\) for any \(p>0\).
Hence, the right hand side of (\ref{euler-lagrange-psi-coordinates}) is bounded in \(L^p\) with \(p>2\).
Again, by the Sobolev embedding theorem we find that \(\psi\in C^{0,\gamma}\) for some \(\gamma>0\).
By elliptic estimates applied to (\ref{euler-lagrange-u-coordinates}) we deduce that \(\phi\in W^{2,p}\)
for any \(p>2\) and thus \(\phi\in C^{1,\alpha}\) for \(0<\alpha<1\). Using elliptic estimates for (\ref{euler-lagrange-psi-coordinates})
we may follow that \(\psi\in C^{1,\gamma}\). By a standard bootstrap argument we then get that \((\phi,\psi)\) is smooth.
\end{proof}
We may summarize our considerations:
\begin{Cor}
Let \((M,h)\) be a closed Riemannian spin surface. Suppose that the pair \((\phi,\psi)\in\chi(M,N)\) is a weak
magnetic Dirac-harmonic map from \(M\) to a compact Riemannian manifold \(N\). Then \(\phi\) is continuous
and thus the pair \((\phi,\psi)\) is smooth.
\end{Cor}

\begin{Bem}
By adjusting the regularity theory for Dirac-harmonic maps derived in \cite{MR2544729}
one may also establish the regularity of magnetic Dirac-harmonic maps.
\end{Bem}

\subsection{Removable Singularity Theorem}\mbox{}\\
In this section we want to establish a \emph{removable singularity theorem} for solutions
of (\ref{euler-lagrange-phi}) and (\ref{euler-lagrange-psi}). It is well known that such a theorem holds
in the limiting cases \(\psi=0\) \cite{MR744314}, \(B=0\) \cite{MR2262709} and \(\psi=0, B=0\) \cite{MR604040}.
Our aim is to show that the theorem is still true for both \(\psi\neq 0\) and \(B\neq 0\).
To this end, we adapt the methods provided for Dirac-harmonic maps in \cite{MR2262709}
to the framework of magnetic Dirac-harmonic maps.

Let us define the following ``energy'', which is the crucial quantity in the context of
removable singularities.
\begin{Dfn}
Let \(U\) be a domain on \(M\). We define the energy of the pair \((\phi,\psi)\) on \(U\)
by
\begin{equation}
E(\phi,\psi,U):=\int_U(|d\phi|^2+|\psi|^4).
\end{equation}
\end{Dfn}
This energy is conformally invariant and thus plays an important role. 
As a next step we want to analyze the local behaviour of magnetic Dirac-harmonic maps
similar to Theorem 4.3 in \cite{MR2262709}. Note that on the unit disc \(D\) the usual Dirac operator satisfies the equation \(\p^2=-\Delta\).
\begin{Satz}
\label{theorem-energy-local}
Let \((M,h)\) be a closed Riemannian spin surface and \((N,g)\) a
compact Riemannian manifold. Assume that the pair \((\phi,\psi)\) is a magnetic Dirac-harmonic map.
There is a small constant \(\epsilon>0\) such that if the pair \((\phi,\psi)\) satisfies
\begin{equation}
\int_D(|d\phi|^2+|\psi|^4)<\epsilon,
\end{equation}
then
\begin{equation}
|d\phi|_{C^k(D_\frac{1}{2})}+|\psi|_{C^k(D_\frac{1}{2})} \leq C(|d\phi|_{L^2(D)}+|\psi|_{L^4(D)}),
\end{equation}
where the constant \(C\) depends on \(|Z|_{L^\infty}, N\) and \(k\).
\end{Satz}
\begin{proof}
Suppose that \(D_\frac{1}{2}\subset D^2\subset D^1\subset D\).
We choose a cut-off function \(\eta:0\leq\eta\leq 1\) with \(\eta|_{D^1}=1\) and \(\operatorname{supp}\eta\subset D\).
By equation (\ref{euler-lagrange-u}) we have
\begin{align*}
|\Delta(\eta\phi)|\leq& c_1(|\phi|+|d\phi|)+(|\sff|_{L^\infty}+|Z|_{L^\infty})|d\phi||d(\eta\phi)|+|\phi d\eta|\\
&\hspace{0.2cm}+c_3\big|\eta|d\phi||\psi|^2\big| \\
\leq& c_1(|\phi|+|d\phi|)+ c_2|d\phi||d(\eta\phi)|+c_3\big|\eta|d\phi||\psi|^2\big|,
\end{align*}
where we set \(c_2:=|\sff|_{L^\infty}+|Z|_{L^\infty}\).
Thus, we find
\[
|\Delta(\eta\phi)|_{L^{\frac{4}{3}}}\leq c_2\big||d\phi||d(\eta\phi)|\big|_{L^{\frac{4}{3}}}+c_1|\phi|_{W^{1,\frac{4}{3}}}+c_3\big|\eta|d\phi||\psi|^2\big|_{L^{\frac{4}{3}}}.
\]
Without loss of generality we assume \(\int_D\phi=0\) such that \(|\phi|_{W^{1,p}}\leq C|d\phi|_{L^p}\) for any \(p>0\).
Moreover, we have by Hölder's inequality
\[
\big||d\phi||d(\eta\phi)|\big|_{L^{\frac{4}{3}}}\leq |\eta\phi|_{W^{1,4}}|d\phi|_{L^2}
\]
such that we may conclude 
\[
|\eta\phi|_{W^{2,\frac{4}{3}}}\leq C(c_2|\eta\phi|_{W^{1,4}}|d\phi|_{L^2}+|d\phi|_{L^{\frac{4}{3}}}+\big|\eta|d\phi||\psi|^2\big|_{L^{\frac{4}{3}}}).
\]
By the Sobolev embedding theorem we find \(|\eta\phi|_{W^{1,4}}\leq c_4|\eta\phi|_{W^{2,\frac{4}{3}}}\) and we may follow
\begin{equation}
\label{estimate-etaphi}
(c_4^{-1}-Cc_2|d\phi|_{L^2})|\eta\phi|_{W^{1,4}}\leq C(|d\phi|_{L^\frac{4}{3}}+\big|\eta|d\phi||\psi|^2\big|_{L^{\frac{4}{3}}}).
\end{equation}
Regarding the last term in (\ref{estimate-etaphi}) we note that
\[
\big|\eta|d\phi||\psi|^2\big|_{L^{\frac{4}{3}}}\leq|\psi|^2_{L^4}|\eta\phi|_{W^{1,4}}+|\psi|^2_{L^4}.
\]
Applying this estimate and choosing \(\epsilon\) small enough, (\ref{estimate-etaphi}) gives
\[
|\eta\phi|_{W^{1,4}}\leq C( |d\phi|_{L^{\frac{4}{3}}}+\sqrt{\epsilon}|\eta\phi|_{W^{1,4}}+|\psi|^2_{L^4}),
\]
which can be rearranged as
\[
|\eta\phi|_{W^{1,4}}\leq C(|d\phi|_{L^\frac{4}{3}}+|\psi|^2_{L^4})\leq \sqrt{\epsilon} C.
\]
Finally, we may deduce that for some \(\epsilon>0\) using the properties of \(\eta\)
\begin{equation}
\label{estimate-step1}
|d\phi|_{L^4(D^1)}\leq C(D^1)(|d\phi|_{L^2(D)}+|\psi|^2_{L^4(D)})\qquad \forall D^1\subset D
\end{equation}
holds.

To establish an estimate for the spinor \(\psi\), we again choose a cut-off function \(\eta\) 
satisfying \(0\leq\eta\leq 1\) with \(\eta|_{D^2}=1\) and \(\operatorname{supp}\eta\subset D\).
Consider the spinor \(\xi:=\eta\psi\) and, using \(\eqref{euler-lagrange-psirq}\), we calculate
\[
\p(\eta\psi)=\eta\p\psi+\nabla\eta\cdot\psi 
=\eta\sff(d\phi(e_\alpha),e_\alpha\cdot\psi)+\nabla\eta\cdot\psi.
\]
Hence, we have
\[
|\xi|_{W^{1,q}(D)}\leq C(\big||d\phi||\eta\psi|\big|_{L^q(D)}+|\psi|_{L^q(D)})
\]
and by Hölder's inequality we can estimate
\[
\big||d\phi||\eta\psi|\big|_{L^q(D)}\leq|d\phi|_{L^2(D)}|\eta\psi|_{L^{q^*}(D)}
\]
with the conjugate Sobolev index \(q^*=\frac{2q}{2-q}\).
By the Sobolev embedding theorem we may then follow
\begin{equation}
\label{nablapsi-l2}
|\xi|_{L^{q^*}(D)}\leq C(|d\phi|_{L^2(D)}|\xi|_{L^{q^*}(D)}+|\psi|_{L^q(D)}). 
\end{equation}
Thus, if the energy \(|d\phi|_{L^2(D)}\) is small enough, we have
\[
|\xi|_{L^{q^*}(D)}\leq C|\psi|_{L^q(D)}.
\]
At this point for any \(p>1\) one can always find some \(q<2\) such that \(p=q^*\),
giving
\begin{equation}
\label{estimate-step2}
|\psi|_{L^q(D^2)}\leq C(D^2)|\psi|_{L^4(D)} \qquad \forall q>1,\qquad D^2\subset D^1\subset D,
\end{equation}
where we again assumed that \(\epsilon\) is small.

Combining (\ref{estimate-step1}), (\ref{estimate-step2}) and for \(\epsilon\) small enough,
we get 
\begin{equation}
|d\phi|_{L^4(D^2)}\leq C(D^2)|d\phi|_{L^2(D)}\qquad \forall D^2\subset D^1\subset D.
\end{equation}
The assertion then follows from an iteration of the procedure presented above,
for more details see the proof of Theorem 3.1 in \cite{MR2176464}.
\end{proof}

A phenomenon connected to harmonic maps from surfaces is the so-called \emph{bubbling}.
Such a bubble arises if one performs a blowup analysis of a sequence of harmonic maps
and we have something similar for magnetic Dirac-harmonic maps.
More precisely, a bubble is a solution of (\ref{euler-lagrange-phi}) and (\ref{euler-lagrange-psi})
with finite energy. Due to the conformal invariance of magnetic Dirac-harmonic maps
such a solution can be interpreted as a magnetic Dirac-harmonic map from \(S^2\setminus\{p\}\to N\)
with finite energy. As a next step we prove how such a singularity can be removed by
exploiting the conformal invariance of our problem.

The behaviour of the pair \((\phi,\psi)\) near a singularity can be described by the following 
\begin{Cor}
\label{corollary-energy-local}
There is an \(\epsilon>0\) small enough such that if the pair \((\phi,\psi)\) is a smooth solution of (\ref{euler-lagrange-phi})
and (\ref{euler-lagrange-psi}) on \(D\setminus\{0\}\) with finite energy \(E(\phi,\psi,D)<\epsilon\), 
then for any \(x\in D_{\frac{1}{2}}\) we have
\begin{eqnarray}
|d\phi(x)||x|&\leq& C\left(\int_{D(2|x|)}|d\phi|^2\right)^\frac{1}{2}, \\
|\psi(x)||x|^\frac{1}{2}+|\nabla\psi(x)||x|^\frac{3}{2}&\leq& C\left(\int_{D(2|x|)}|\psi|^4\right)^\frac{1}{4}.
\end{eqnarray}
\end{Cor}
\begin{proof}
Fix any \(x_0\in D\setminus\{0\}\) and define \((\tilde{\phi},\tilde{\psi})\) by
\[
\tilde{\phi}(x):=\phi(x_0+|x_0|x) \textrm{ and } \tilde{\psi}(x):=|x_0|^\frac{1}{2}\psi(x_0+|x_0|x).
\]
It is easy to see that \((\tilde{\phi},\tilde{\psi})\) is a smooth solution of (\ref{euler-lagrange-phi}) and (\ref{euler-lagrange-psi}) on \(D\)
with \(E(\tilde{\phi},\tilde{\psi},D)<\epsilon\). By application of Theorem \ref{theorem-energy-local}, we have
\begin{eqnarray*}
|d\tilde{\phi}|_{L^\infty(D_\frac{1}{2})}&\leq& C|d\tilde{\phi}|_{L^2(D)}, \\
|\tilde{\psi}|_{C^1(D_\frac{1}{2})}&\leq& C|\tilde{\psi}|_{L^4(D)}
\end{eqnarray*}
and scaling back yields the assertion.
\end{proof}

\begin{Lem}
\label{lemma-polar-coordinates}
Let \((\phi,\psi)\) be a smooth solution of (\ref{euler-lagrange-phi}) and (\ref{euler-lagrange-psi}) on \(D\setminus\{0\}\) with
\(E(\phi,\psi,D)<\epsilon\). Then, we have
\begin{eqnarray}
\int_0^{2\pi}\frac{1}{r^2}\big|\frac{\partial\phi}{\partial\theta}\big|^2d\theta
&=&\int_0^{2\pi}\big|\frac{\partial\phi}{\partial r}\big|^2d\theta+
\nonumber\int_0^{2\pi}\langle\psi,\partial_r\cdot\frac{\tilde{\nabla}\psi}{\partial r}\rangle d\theta\\
&=&\int_0^{2\pi}\big|\frac{\partial\phi}{\partial r}\big|^2d\theta-
\int_0^{2\pi}\frac{1}{r^2}\langle\psi,\partial_\theta\cdot\frac{\tilde{\nabla}\psi}{\partial\theta}\rangle d\theta,
\end{eqnarray}
where \((r,\theta)\) are polar coordinates in \(D\) centered at \(0\).
\end{Lem}
\begin{proof}
Applying the previous Corollary, it can be shown that \(|T(z)|\leq Cz^{-2}\), where \(T(z)\) is the Hopf differential
defined in (\ref{hopf-differential}). Moreover, it follows that \(\int_D|T(z)|<\infty\), hence \(zT(z)\) is a holomorphic function on \(D\). 
The claim then follows from Cauchy's integral theorem 
\[
0=\int_{|z|=r}\operatorname{Im}(zT(z))dz=\int_{0}^{2\pi}\operatorname{Re}(z^2T(z))d\theta
\]
and changing to polar coordinates. 
For more details see the proof of Lemma 4.5 in \cite{MR2262709}.
\end{proof}
Additionally, by integration the last Lemma also gives
\begin{equation}
\int_D\big|\frac{\partial\phi}{\partial r}\big|^2-\int_D\frac{1}{r^2}\big|\frac{\partial\phi}{\partial\theta}\big|^2
=-\int_D\langle\psi,\partial_r\cdot\frac{\tilde{\nabla}\psi}{\partial r}\rangle:=I.
\end{equation}
For the Dirichlet energy in polar coordinates, we have
\[
E_r(\phi):=\int_{|z|=r}|d\phi|^2=\int_{|z|=r}\big|\frac{\partial\phi}{\partial r}\big|^2+\int_{|z|=r}\frac{1}{r^2}\big|\frac{\partial\phi}{\partial\theta}\big|^2.
\]
Set \(I_r(\psi)=-\int_{|z|=r}\langle\psi,\partial_r\cdot\frac{\tilde{\nabla}\psi}{\partial r}\rangle\), then we get
\begin{eqnarray}
\label{proof-singularity-E}\int_{|z|=r}\big|\frac{\partial\phi}{\partial r}\big|^2&=&\frac{1}{2}E_r(\phi)+\frac{1}{2}I_r(\psi), \\
\label{proof-singularity-I}\int_{|z|=r}\frac{1}{r^2}\big|\frac{\partial\phi}{\partial\theta}\big|^2&=&\frac{1}{2}E_r(\phi)-\frac{1}{2}I_r(\psi).
\end{eqnarray}
Before we turn to the removable singularity theorem, let us state an auxiliary 
\begin{Lem}[Elliptic estimate with boundary]
\label{estimate-boundary}
Suppose that \(\psi\) is a solution of
\begin{eqnarray}
\p\psi&=&f \qquad \text{ on } D, \\
\nonumber\psi&=&g \qquad \text{ on } \partial D.
\end{eqnarray}
Assume that \(f\in L^p(D,\Sigma M)\) and \(g\in W^{1,p}(\partial D,\Sigma M)\) for some \(p>1\), then
the following estimate holds
\begin{equation}
|\psi|_{W^{1,p}(D)}\leq C(|g|_{W^{1,p}(\partial D)}+|f|_{L^p(D)}),
\end{equation}
with the positive constant \(C=C(p)\).
\end{Lem}
A proof can be found in \cite{MR2390834}, Lemma 2.2.

Now we are in the position to state the
\begin{Satz}[Removable Singularity Theorem]
Let \((\phi,\psi)\) be a solution of (\ref{euler-lagrange-phi}) and (\ref{euler-lagrange-psi}), which is smooth on \(U\setminus\{p\}\) for some \(p\in U\).
If \((\phi,\psi)\) has finite energy, then \((\phi,\psi)\) extends to a smooth solution on \(U\). 
\end{Satz}

\begin{proof}
By rescaling we may assume that 
\[
\int_{D(2)}(|d\phi|^2+|\psi|^4)<\epsilon,
\]
where \(\epsilon\) is given by Theorem \ref{theorem-energy-local}.
We will approximate \(\phi\) by a function \(q=q(r)\) that depends only on the radial coordinate
and which is piecewise linear in \(\log r\). For \(m\geq 1\) set
\[
q(2^{-m})=\frac{1}{2\pi}\int_0^{2\pi}\phi(2^{-m},\theta)d\theta.
\]
Then \(q\) is harmonic for \(r\in(2^{-m},2^{-m+1}),m\geq 1\), see for example \cite{MR2431658}, p.130.
Moreover, we have 
\begin{equation}
\label{proof-singularity-a}
\int_D|dq-d\phi|^2=\int_{r=1}\langle q-\phi,\frac{\partial\phi}{\partial r}\rangle-\int_D\langle q-\phi,\Delta(q-\phi)\rangle
\end{equation}
and
\begin{equation}
|q-\phi|_{L^\infty}:=|q-\phi|_{C^0(D)}\leq 2^3\sqrt{\epsilon}.
\end{equation}
Using (\ref{euler-lagrange-u}), we may then estimate
\begin{eqnarray}
\label{proof-singularity-b}
\nonumber \left|\int_{D}\langle\Delta \phi,\phi-q\rangle\right|&\leq& (|\sff|_{L^\infty}+|Z|_{L^\infty})|\phi-q|_{L^\infty(D)}\int_D|d\phi|^2\\
\nonumber &&+C|\phi-q|_{L^\infty(D)}\int_D|\psi|^2|d\phi|\\
 &\leq&\delta\int_D|d\phi|^2+C\sqrt{\epsilon}\int_D|\psi|^2|d\phi|.
\end{eqnarray}
with \(2^3\sqrt{\epsilon}(|\sff_{L^\infty}|+|Z|_{L^\infty})\leq\delta\).
We may estimate the first term on the right hand side of (\ref{proof-singularity-a}) as
\begin{equation}
\label{proof-singularity-c}
\int_{r=1}\langle q-\phi,\frac{\partial\phi}{\partial r}\rangle\leq\bigg(\int_{r=1}\big|\frac{\partial\phi}{\partial\theta}\big|^2\bigg)^\frac{1}{2}\bigg(\int_{r=1}\big|\frac{\partial\phi}{\partial r}\big|^2\bigg)^\frac{1}{2}.
\end{equation}
By Lemma \ref{lemma-polar-coordinates} and the fact that the function \(q\) does not depend on \(\theta\), we find
\begin{equation}
\label{proof-singularity-d}
\frac{1}{2}\int_{D_1}|d\phi|^2-\frac{1}{2}I=\int_{D_1}\frac{1}{r^2}\big|\frac{\partial\phi}{\partial\theta}\big|^2\leq \int_{D_1}|d(\phi-q)|^2.	 
\end{equation}
Inserting (\ref{proof-singularity-b}), (\ref{proof-singularity-c}) and (\ref{proof-singularity-d}) into (\ref{proof-singularity-a}) then yields 
\begin{eqnarray}
\frac{1}{2}\int_{D_1}|d\phi|^2-\frac{1}{2}I&\leq&\bigg(\int_{r=1}\big|\frac{\partial\phi}{\partial\theta}\big|^2\bigg)^\frac{1}{2}\bigg(\int_{r=1}\big|\frac{\partial\phi}{\partial r}\big|^2\bigg)^\frac{1}{2} \\
\nonumber&&+\delta\int_D|d\phi|^2+C\sqrt{\epsilon}\int_D|\psi|^2|d\phi|.
\end{eqnarray}
By (\ref{proof-singularity-E}) and (\ref{proof-singularity-I}) we find
\begin{eqnarray*}
\bigg(\int_{r=1}\big|\frac{\partial\phi}{\partial\theta}\big|^2\bigg)^\frac{1}{2}\bigg(\int_{r=1}\big|\frac{\partial\phi}{\partial r}\big|^2\bigg)^\frac{1}{2}
&=&\frac{1}{2}(E_1-I_1)^\frac{1}{2}(E_1+I_1)^\frac{1}{2}\\
&\leq&\frac{1}{2}E_1=\frac{1}{2}\int_{r=1}|d\phi|^2
\end{eqnarray*}
and it follows that
\begin{equation}
(1-2\delta)\int_D|d\phi|^2\leq\int_{r=1}|d\phi|^2+2C\sqrt{\epsilon}\int_D|\psi|^2|d\phi|-\int_D\langle\psi,\partial_r\cdot\frac{\nabla}{\partial r}\psi\rangle.
\end{equation}
By rescaling and some manipulations on the right hand side, we arrive at
\begin{eqnarray}
\label{F-a}
(1-2\delta)\int_{D_r}|d\phi|^2&\leq& r\int_{\partial D_r}|d\phi|^2 \\
\nonumber &&+C\big(\sqrt{\epsilon}\int_{D_r}|d\phi|^2+\int_{D_r}|\psi|^4+\int_{D_r}|\nabla\psi|^\frac{4}{3}\big).
\end{eqnarray}
Our next aim is to estimate all contributions on the right hand side in terms of integrals over the boundary \(\partial D_r\).
To this end, we apply Lemma \ref{estimate-boundary} to the Euler-Lagrange equation (\ref{euler-lagrange-psirq}) for \(\psi\).
Thus, we find
\begin{eqnarray*}
|\psi|_{W^{1,\frac{4}{3}}(D)}&\leq& C\big(\big||d\phi||\psi|\big|_{L^\frac{4}{3}(D)}+|\psi|_{W^{1,\frac{4}{3}}(\partial D)}\big) \\
&\leq&C\big(|d\phi|_{L^2(D)}|\psi|_{L^4(D)}+|\psi|_{W^{1,\frac{4}{3}}(\partial D)}\big).
\end{eqnarray*}
By the Sobolev embedding theorem and the smallness of \(|d\phi|_{L^2(D)}\) we can achieve that
\[
|\psi|_{L^4(D)}\leq C\big(|\nabla\psi|_{L^\frac{4}{3}(\partial D)}+|\psi|_{L^4(\partial D)}\big).
\]
By rescaling, we have for \(0<r\leq 1\) 
\begin{equation}
\label{F-b}
\int_{D_r}|\psi|^4\leq Cr\int_{\partial D_r}|\nabla\psi|^\frac{4}{3}+Cr\int_{\partial D_r}|\psi|^4.
\end{equation}
By the same methods one can derive
\begin{equation*}
|\nabla\psi|_{L^\frac{4}{3}(D)}\leq C\big(|d\phi|_{L^2(D)}|\psi|_{L^4(D)}+|\psi|_{W^{1,\frac{4}{3}}(\partial D)}\big)
\end{equation*}
and again by the smallness of \(|d\phi|_{L^2(D)}\) and a scaling argument we find
\begin{equation}
\label{F-c}
\int_{D_r}|\nabla\psi|^\frac{4}{3}\leq C\sqrt{\epsilon}\int_{D_r}|\psi|^4+Cr\int_{\partial D_r}|\nabla\psi|^\frac{4}{3}+Cr\int_{\partial D_r}|\psi|^4.
\end{equation}
Combining (\ref{F-a}), (\ref{F-b}) and (\ref{F-c}), we have for any \(0<r\leq 1\) and some constant \(C>0\) 
\begin{eqnarray}
&&\int_{D_r}|d\phi|^2+\int_{D_r}|\psi|^4+\int_{D_r}|\nabla\psi|^\frac{4}{3} \\
\nonumber&&\hspace{2cm}\leq Cr\big(\int_{\partial D_r}|d\phi|^2+\int_{\partial D_r}|\psi|^4+\int_{\partial D_r}|\nabla\psi|^\frac{4}{3}\big).
\end{eqnarray}
Setting \(F(r):=\int_{D_r}|d\phi|^2+\int_{D_r}|\psi|^4+\int_{D_r}|\nabla\psi|^\frac{4}{3}\), we obtain
\begin{equation}
F(r)\leq Cr F'(r), 
\end{equation}
which can be integrated as
\begin{equation}
F(r)\leq F(1)r^\frac{1}{C}.
\end{equation}
Hence, there are some \(a>1\) and \(2>b>\frac{4}{3}\) such that
\begin{equation}
\phi\in W^{1,2a}(D,\R^q),\qquad \psi\in W^{1,b}(D,\Sigma M\otimes T\R^q) 
\end{equation}
holds. By applying an iterated bootstrap procedure one can now show that the pair \((\phi,\psi)\) is smooth on \(D\).
For more details, see the proof of Theorem 4.6 in \cite{MR2262709}.
\end{proof}
\bibliographystyle{amsplain}
\bibliography{mybib}

\def\cprime{$'$}
\providecommand{\bysame}{\leavevmode\hbox to3em{\hrulefill}\thinspace}
\providecommand{\MR}{\relax\ifhmode\unskip\space\fi MR }
\providecommand{\MRhref}[2]{%
  \href{http://www.ams.org/mathscinet-getitem?mr=#1}{#2}
}
\providecommand{\href}[2]{#2}
\begin{thebibliography}{10}

\bibitem{MR1910266}
Orlando Alvarez and I.~M. Singer, \emph{Beyond the elliptic genus}, Nuclear
  Phys. B \textbf{633} (2002), no.~3, 309--344.

\bibitem{MR3070562}
Bernd Ammann and Nicolas Ginoux, \emph{Dirac-harmonic maps from index theory},
  Calc. Var. Partial Differential Equations \textbf{47} (2013), no.~3-4,
  739--762.

\bibitem{MR1168525}
Fabrice Bethuel, \emph{Un r\'esultat de r\'egularit\'e pour les solutions de
  l'\'equation de surfaces \`a courbure moyenne prescrite}, C. R. Acad. Sci.
  Paris S\'er. I Math. \textbf{314} (1992), no.~13, 1003--1007.

\bibitem{phd}
Volker Branding, \emph{The evolution equations for {D}irac-harmonic maps,
  {P}h{D} {T}hesis}, Universitaet Potsdam (2013).

\bibitem{MR2370260}
Q.~Chen, J.~Jost, and G.~Wang, \emph{Liouville theorems for {D}irac-harmonic
  maps}, J. Math. Phys. \textbf{48} (2007), no.~11, 113517, 13.

\bibitem{BVP}
Q.~Chen, J.~Jost, G.~Wang, and M.~Zhu, \emph{The boundary value problem for
  dirac-harmonic maps}, J. of the EMS \textbf{15} (2013), no.~3, 997--1031.

\bibitem{MR2176464}
Qun Chen, J{\"u}rgen Jost, Jiayu Li, and Guofang Wang, \emph{Regularity
  theorems and energy identities for {D}irac-harmonic maps}, Math. Z.
  \textbf{251} (2005), no.~1, 61--84.

\bibitem{MR2262709}
\bysame, \emph{Dirac-harmonic maps}, Math. Z. \textbf{254} (2006), no.~2,
  409--432.

\bibitem{MR2390834}
Qun Chen, J{\"u}rgen Jost, and Guofang Wang, \emph{Nonlinear {D}irac equations
  on {R}iemann surfaces}, Ann. Global Anal. Geom. \textbf{33} (2008), no.~3,
  253--270.

\bibitem{MR3044133}
\bysame, \emph{The maximum principle and the {D}irichlet problem for
  {D}irac-harmonic maps}, Calc. Var. Partial Differential Equations \textbf{47}
  (2013), no.~1-2, 87--116.

\bibitem{MR1331835}
Philippe Chon{\'e}, \emph{A regularity result for critical points of
  conformally invariant functionals}, Potential Anal. \textbf{4} (1995), no.~3,
  269--296.

\bibitem{MR2648325}
J{\"u}rgen Fuchs, Thomas Nikolaus, Christoph Schweigert, and Konrad Waldorf,
  \emph{Bundle gerbes and surface holonomy}, European {C}ongress of
  {M}athematics, Eur. Math. Soc., Z\"urich, 2010, pp.~167--195.

\bibitem{MR636440}
Michael Gr{\"u}ter, \emph{Regularity of weak {$H$}-surfaces}, J. Reine Angew.
  Math. \textbf{329} (1981), 1--15.

\bibitem{MR744314}
\bysame, \emph{Conformally invariant variational integrals and the removability
  of isolated singularities}, Manuscripta Math. \textbf{47} (1984), no.~1-3,
  85--104.

\bibitem{MR1913803}
Fr{\'e}d{\'e}ric H{\'e}lein, \emph{Harmonic maps, conservation laws and moving
  frames}, second ed., Cambridge Tracts in Mathematics, vol. 150, Cambridge
  University Press, Cambridge, 2002, Translated from the 1996 French original,
  With a foreword by James Eells.

\bibitem{MR1242500}
C.~M. Hull, G.~Papadopoulos, and P.~K. Townsend, \emph{Potentials for {$(p,0)$}
  and {$(1,1)$} supersymmetric sigma models with torsion}, Phys. Lett. B
  \textbf{316} (1993), no.~2-3, 291--297.

\bibitem{kohdiss}
Dennis Koh, \emph{The evolution equation for closed magnetic geodesics},
  Dissertation, Universitätsverlag Potsdam (2008).

\bibitem{MR1031992}
H.~Blaine Lawson, Jr. and Marie-Louise Michelsohn, \emph{Spin geometry},
  Princeton Mathematical Series, vol.~38, Princeton University Press,
  Princeton, NJ, 1989.

\bibitem{MR2431658}
Fanghua Lin and Changyou Wang, \emph{The analysis of harmonic maps and their
  heat flows}, World Scientific Publishing Co. Pte. Ltd., Hackensack, NJ, 2008.

\bibitem{MR2285745}
Tristan Rivi{\`e}re, \emph{Conservation laws for conformally invariant
  variational problems}, Invent. Math. \textbf{168} (2007), no.~1, 1--22.

\bibitem{MR604040}
J.~Sacks and K.~Uhlenbeck, \emph{The existence of minimal immersions of
  {$2$}-spheres}, Ann. of Math. (2) \textbf{113} (1981), no.~1, 1--24.

\bibitem{MR2788659}
Matthias Schneider, \emph{Closed magnetic geodesics on {$S^2$}}, J.
  Differential Geom. \textbf{87} (2011), no.~2, 343--388.

\bibitem{1306.4260}
Benjamin Sharp and Miaomiao Zhu, \emph{Regularity at the free boundary for
  {D}irac-harmonic maps from surfaces}, arXiv:1306.4260 (2013).

\bibitem{MR2681604}
Konrad Waldorf, \emph{Surface holonomy}, Handbook of pseudo-{R}iemannian
  geometry and supersymmetry, IRMA Lect. Math. Theor. Phys., vol.~16, Eur.
  Math. Soc., Z\"urich, 2010, pp.~653--682.

\bibitem{MR2544729}
Changyou Wang and Deliang Xu, \emph{Regularity of {D}irac-harmonic maps}, Int.
  Math. Res. Not. IMRN (2009), no.~20, 3759--3792.

\bibitem{MR3018163}
Deliang Xu and Zhengxiang Chen, \emph{Regularity for {D}irac-harmonic map with
  {R}icci type spinor potential}, Calc. Var. Partial Differential Equations
  \textbf{46} (2013), no.~3-4, 571--590.

\bibitem{MR2267756}
Liang Zhao, \emph{Energy identities for {D}irac-harmonic maps}, Calc. Var.
  Partial Differential Equations \textbf{28} (2007), no.~1, 121--138.

\bibitem{MR2506243}
Miaomiao Zhu, \emph{Regularity for weakly {D}irac-harmonic maps to
  hypersurfaces}, Ann. Global Anal. Geom. \textbf{35} (2009), no.~4, 405--412.

\end{thebibliography}

\end{document}